\documentclass[11pt]{amsart}
\usepackage{amsmath,amsthm,amsfonts,amssymb,latexsym,epsfig}

\newtheorem{theorem}{Theorem}[section]

\newtheorem{lemma}{Lemma}[section]

\newtheorem{cor}{Corollary}[section]

\numberwithin{equation}{section}

\theoremstyle{definition}

\theoremstyle{remark}

\begin{document}
\title{A Note on Sums of Powers}
\author{Peng Gao}
\address{Department of Computer and Mathematical Sciences,
University of Toronto at Scarborough, 1265 Military Trail, Toronto
Ontario, Canada M1C 1A4} \email{penggao@utsc.utoronto.ca}
%%\date{\today}
\date{April 17, 2007.}

%%-------------------------------------------------------------------------------------------

\begin{abstract}
 We improve a result of Bennett concerning 
 certain sequences involving sums of powers of positive integers.
\end{abstract}

\maketitle
%%-------------------------------------------------------------------------------------------
\section{Introduction}
\label{sec 1} \setcounter{equation}{0}
%%-------------------------------------------------------------------------------------------
   Estimations of sums of powers of positive integers have important applications in
  the study of $l^p$ norms of weighted mean matrices, we leave
  interested readers the recent papers \cite{G} and \cite{Be1}
  for more details in this direction. There are many inequalities for sequences involving sums of powers of
   positive integers in the literature and we shall also refer the interested readers to the recent papers \cite{Be1}, \cite{Be2}, 
\cite{G1} as well as the references therein for some results in this area.

   In this note, we are interested in certain inequalities involving the following sequence: 
$\{ P_n(r) | n=1, 2, 3, \ldots \}$, where $r$ is any real number and 
\begin{equation*}
  P_{n}(r)= \left(\frac {1}{n} \sum_{i=1}^{n}i^r\bigg/
\frac {1}{n+1}\sum_{i=1}^{n+1}i^r\right)^{1/r}, \hspace{0.1in} r
\neq 0; \hspace{0.1in} P_{n}(0)=\frac
   {\sqrt[n]{n!}}{\sqrt[n+1]{(n+1)!}}.
\end{equation*}
  We note that for $r > 0$,
  the following inequalities hold:
\begin{equation}
\label{1}
   \frac {n}{n+1}= \lim_{r \rightarrow +\infty}P_n(r) < P_n(r) < P_n(0).
\end{equation}
  The left-hand side inequality above is known as Alzer's inequality
\cite{alz}, and the right-hand side inequality above is known as
Martins' inequality \cite{Mar}.  Alzer also considered inequalities
satisfied by $P_n(r)$ for $r<0$ in \cite{alz1} and he showed
\cite[Theorem 2.3]{alz1}:
\begin{equation}
\label{2}
   P_n(0) \leq P_n(r) \leq \lim_{r \rightarrow -\infty}P_n(r)=1.
\end{equation}

  Bennett \cite{Be} proved that for $r \geq 1$,
\begin{equation}
\label{1.5}
  P_n(r) \leq P_n(1)=\frac
   {n+1}{n+2}
\end{equation}
  with the above inequality reversed when $0< r \leq 1$. This inequality
  and inequalities \eqref{1}-\eqref{2} suggest that $P_n(r)$ is a decreasing
  function of $r$. Recently, Bennett \cite{Be2} proved this for $r \leq
  1$ and the author gave another proof in \cite{G1}. Bennett further asked, using his notation
  in \cite{Be2}, to decide whether the sequence 
$(1^r, 2^r, 3^r, \ldots)$ is meaningful for any $r > 1$ or not (\cite[Problem 1]{Be2}), which is equivalent to asking for 
whether $P_n(r)$ is a decreasing
  function of $r$ for any $r > 1$ or not.  
  It is our goal in this note to give a weaker result related to Bennett's question above by proving 
  the following:
\begin{theorem}
\label{thm4}
  The sequence
\begin{equation*}
   \frac {\Big ( \sum^n_{i=1}i^r \Big
  )^{\alpha}}{\sum^n_{i=1}i^{\alpha(r+1)-1}}, \hspace {0.1in} n=1,
  2, 3, \ldots,
\end{equation*}
   is decreasing for $r \geq 1, \alpha \geq 2$.
\end{theorem}  

  We note here that Theorem \ref{thm4} improves a result of Bennett \cite[Theorem 12]{Be1},
  which established the case $\alpha=2$ of Theorem \ref{thm4}. We also note that one can readily 
  deduce from Theorem \ref{thm4} using an argument similar to the discussion in the paragraph below Corollary 3.1 
in \cite{G1} the following 
\begin{cor}
\label{cor10}
  For any fixed integer $n \geq 1$, $P_n(r) \geq P_n(r')$ for $r' \geq 2r+1, r \geq
   1$.
\end{cor}

%%-------------------------------------------------------------------------------------------
\section{Lemmas}
\label{sec 2} \setcounter{equation}{0}
%%-------------------------------------------------------------------------------------------
\begin{lemma}[{\cite[Lemma 2.1]{Xu}}]
\label{lem3}
   Let $\{B_n \}^{\infty}_{n=1}$ and $\{C_n \}^{\infty}_{n=1}$ be strictly increasing positive sequences with
   $B_1/B_2 \leq C_1 / C_2$. If for any integer $n \geq 1$,
\begin{equation*}
  \frac {B_{n+1}-B_n}{B_{n+2}-B_{n+1}} \leq  \frac
  {C_{n+1}-C_n}{C_{n+2}-C_{n+1}}.
\end{equation*}
  Then $B_{n}/B_{n+1} \leq C_{n} / C_{n+1}$ for any integer $n \geq 1$.
\end{lemma}

\begin{lemma}
\label{lem4}
  For $r \geq 2, x>0, y>0$, let
\begin{equation*}
    D_r (x, y) = \frac {x^r - y^r}{ x-y}, \hspace{0.1in} x \neq y; \hspace{0.1in} D_r(x,x)=rx^{r-1}.
\end{equation*}
   Then for positive numbers $a, b, c, d$ satisfying $a \geq \max ( b, c, d)$ and $a+b \geq c+d$, we have
\begin{equation*}
   D_r(a, b) \geq D_r(c, d).
\end{equation*}
\end{lemma}
\begin{proof}
  We may assume $c \geq d$ here and note that $D_r(x,y)$ is an increasing function of $x$ (or $y$) for fixed $y$ (or $x$).
  It follows from this that if $b \geq d$, then $D_r(a, b) \geq D_r(c, b) \geq D_r(c, d)$. Otherwise by our assumption, 
  one can find a positive number $a'$ such that $a \geq a' \geq \max ( b, c, d)$ and $a'+b = c+d$. 
  
    We now recall from the theory of majorization that for two positive real finite sequences
   ${\bf x}=(x_1, x_2, \ldots, x_n)$ and ${\bf y}=(y_1, y_2,
   \ldots, y_n)$, ${\bf x}$ is said to be
   majorized by ${\bf y}$ if for all convex functions $f$, we have
\begin{equation*}
%%\label{0}
   \sum_{j=1}^{n}f(x_j) \leq \sum_{j=1}^{n}f(y_j).
\end{equation*}

    We write ${\bf x} \leq_{maj} {\bf y}$ if this occurs and the
    majorization principle states that if $(x_j)$ and $(y_j)$ are
    decreasing, then ${\bf x} \leq_{maj} {\bf y}$ is equivalent to
\begin{eqnarray*}
\label{10}
   x_1+x_2+\ldots+x_j & \leq & y_1+y_2+\ldots+y_j ~~(1 \leq j \leq
   n-1),
     \\
    x_1+x_2+\ldots+x_n & = & y_1+y_2+\ldots+y_n ~~(n \geq 0).
\end{eqnarray*}
    We refer the reader to \cite[Sect. 1.30]{B&B} for a simple proof of
    this.
   
    Now let $I \subset (0, +\infty)$ be an open interval and denote
  $I^n = I \times I \times \cdots \times I$ ( $n$ copies). We recall a function $f: I^n
\rightarrow R$ is said to be Schur convex if $f({\bf x}) \leq f({\bf y})$ for any two sequences ${\bf x}, {\bf y} \in I^n$ 
with ${\bf x} \leq_{maj} {\bf y}$.   If $f$
also has continuous partial derivatives on $I^n$, then $f$ is Schur
convex if and only if (see \cite[p. 57]{M&O})
\begin{equation}
\label{8.17}
   (x_i-x_j)(\frac {\partial f}{\partial x_i}-\frac {\partial f}{\partial
x_j}) \geq 0.
\end{equation}
  
  Back to our situation, we apply the notion of majorization to write $(c, d) \leq_{maj} (a', b)$ 
  and we next show that $D_r(x, y)$ satisfies the criterion \eqref{8.17} 
  on $(0, +\infty) \times (0, +\infty)$. For this, we may assume $x > y$ here and then
  it is easy to see that it suffices to show
 \begin{equation*}
   \frac {x^r-y^r}{x-y}=\frac r{x-y}\int^x_yt^{r-1}dt \leq \frac {x^{r-1}+y^{r-1}}{2}.
 \end{equation*}
   The inequality above now follows from Hadamard's inequality which asserts that for a continuous 
  convex function $h(x)$ on an interval $[e, f]$,
\begin{equation*}
%%\label{3.6}
   h(\frac {e+f}2) \leq \frac {1}{f-e}\int^f_e h(x)dx \leq \frac
   {h(e)+h(f)}{2}.
\end{equation*}
  It follows that $D_r(x, y)$ is Schur convex on $(0, +\infty) \times (0, +\infty)$ so that 
  $D_r(a', b) \geq D_r(c, d)$.  As $D_r(a, b) \geq D_r(a', b)$, this completes the proof. 
\end{proof}

\begin{lemma}
\label{lem5}
  For $r \geq 1, \alpha \geq 1$, let
 \begin{equation*}
  g_r(\alpha) =1+2^{\alpha(r+1)-1}-(1+2^r)^{\alpha}.
 \end{equation*}
   Then $g_r(\alpha) \geq 0$ for $r\geq 1, \alpha \geq 2$.
\end{lemma}
\begin{proof}
  We may assume $r \geq 1$ is being fixed and regard $g_r(\alpha)$ as a function of $\alpha$. Then
\begin{equation*}
  g'_r(\alpha)=(\ln 2^{r+1}) 2^{\alpha(r+1)-1} -\ln (1+2^r)(1+2^r)^{\alpha}. 
\end{equation*}
   From this we see that $g'_r(\alpha)=0$ has at most one positive root. Note that $g_r(2) \geq 0$ and 
$\lim_{\alpha \rightarrow +\infty}g_r(\alpha) = +\infty$, it thus suffices to show that $g'_r(2) > 0$.
Note that $g'_r(2)=f(2^r)$, where
\begin{equation*}
  f(x)=\ln (2x) (2x^2) -\ln (1+x)(1+x)^2. 
\end{equation*} 
  As it is easy to check that $f(2) >0, f'(2)>0$, it suffices to show that $f''(x) \geq 0$ for $x \geq 2$. 
  Calculation yields:
\begin{equation*}
  f''(x) =3+4\ln 2+2 \Big(\ln x^2-\ln (1+x) \Big ) > 0.
\end{equation*} 
  The last inequality follows from $x^2 > 1+x$ when $x \geq 2$ and this completes the proof.
\end{proof} 

%%-------------------------------------------------------------------------------------------
\section{Proof of Theorem \ref{thm4}}
\label{sec 3} \setcounter{equation}{0}
%%-------------------------------------------------------------------------------------------
   We need to show that for $n \geq 1$, $r \geq 1, \alpha \geq 2$,
\begin{equation*}
  \frac {\Big ( \sum^n_{i=1}i^r \Big
  )^{\alpha}}{\sum^n_{i=1}i^{\alpha(r+1)-1}} \geq \frac {\Big ( \sum^{n+1}_{i=1}i^r \Big
  )^{\alpha}}{\sum^{n+1}_{i=1}i^{\alpha(r+1)-1}}.
\end{equation*}
  When $n=1$, this
  follows from Lemma \ref{lem5}. Now by Lemma \ref{lem3}, it suffices to show for 
  $n \geq 1$, $r \geq 1, \alpha \geq 2$,
\begin{equation*}
%%\label{3.10}
   \frac
  {\Big ( \sum^{n+1}_{i=1}i^r \Big
  )^{\alpha}-\Big ( \sum^n_{i=1}i^r \Big
  )^{\alpha}}{(n+1)^{\alpha(r+1)-1}} \geq \frac
  {\Big ( \sum^{n+2}_{i=1}i^r \Big
  )^{\alpha}-\Big ( \sum^{n+1}_{i=1}i^r \Big
  )^{\alpha}}{(n+2)^{\alpha(r+1)-1}}.
\end{equation*}
  We can rewrite the above inequality as $D_{\alpha}(a, b) \geq D_{\alpha}(c, d)$, where
 \begin{equation*}
   a= \frac {\sum^{n+1}_{i=1}i^r}{(n+1)^{r+1}}, b= \frac {\sum^{n}_{i=1}i^r}{(n+1)^{r+1}}, 
c=\frac {\sum^{n+2}_{i=1}i^r}{(n+2)^{r+1}}, d=\frac {\sum^{n+1}_{i=1}i^r}{(n+2)^{r+1}}.
 \end{equation*}
 It is easy to see that $a \geq \max (b, d)$ and $a \geq c$ is equivalent to $P_n(r) \geq P_n(0)$, which follows from \eqref{1}.
 Thus our theorem will follow from Lemma \ref{lem4} provided that we show $a+b \geq c+d$ here, which is
 \begin{equation}
 \label{3.11}
     \frac {\sum^{n+1}_{i=1}i^r+\sum^{n}_{i=1}i^r}{(n+1)^{r+1}}  \geq  
  \frac {\sum^{n+2}_{i=1}i^r+\sum^{n+1}_{i=1}i^r}{(n+2)^{r+1}}.
 \end{equation}
   On setting $B_n=n^{r+1}$ and $C_n=\sum^{n}_{i=1}i^r+\sum^{n-1}_{i=1}i^r$ (where we take the empty sum to be $0$) 
 in Lemma \ref{lem3},  it is easy to see that $B_1/B_2 \leq C_1/C_2$. Hence inequality \eqref{3.11} will
 follow from Lemma \ref{lem3} if we can show for $n \geq 1$,
\begin{equation*}
  \frac {(n+1)^{r}+n^{r}}{(n+1)^{r+1}-n^{r+1}}  \geq  
 \frac {(n+2)^{r}+(n+1)^{r}}{(n+2)^{r+1}-(n+1)^{r+1}} .
\end{equation*}
   On setting $x=n/(n+1)$, it is easy to see that one can deduce the above inequality by showing the
  following function is decreasing for $0 < x < 1$:
\begin{equation*}
  f(x)=\frac {(1-x)(1+x^r)}{1-x^{r+1}}.
\end{equation*}
  Calculation yields
\begin{equation*}
  f'(x)=\frac {x^{2r}-rx^{r+1}+rx^{r-1}-1}{(1-x^{r+1})^2}.
\end{equation*}
  It is easy to see that the function $x \mapsto x^{2r}-rx^{r+1}+rx^{r-1}-1$ is an increasing function of $0 < x < 1$ 
  with value $0$ when $x=1$ for any fixed $r \geq 1$.  This implies that $f'(x) \leq 0$ for $0 < x < 1$ and 
  this completes the proof.

%%-----------------------------------------------------------------------------------------

\end{document}